\documentclass[letterpaper, 11pt,  reqno]{amsart}

\usepackage{amsmath,amssymb,amscd,amsthm,amsxtra, esint}
\allowdisplaybreaks[2]

\newtheorem{theorem}{Theorem} [section]

\newtheorem{proposition}[theorem]{Proposition}
\newtheorem{remark}[theorem]{Remark}
\newtheorem{example}{Example}

\newcommand{\noi}{\noindent}
\newcommand{\Z}{\mathbb{Z}}
\newcommand{\R}{\mathbb{R}}
\newcommand{\C}{\mathbb{C}}
\newcommand{\T}{\mathbb{T}}

\newcommand{\al}{\alpha}
\newcommand{\be}{\beta}
\newcommand{\dl}{\delta}

\newcommand{\Dl}{\Delta}
\newcommand{\eps}{\varepsilon}

\newcommand{\g}{\gamma}

\newcommand{\ld}{\lambda}

\newcommand{\s}{\sigma}
\newcommand{\Si}{\Sigma}
\newcommand{\ft}{\widehat}

\newcommand{\wt}{\widetilde}
\newcommand{\cj}{\overline}
\newcommand{\dx}{\partial_x}

\newcommand{\LRA}{\Longrightarrow}

\newcommand{\jb}[1]
{\langle #1 \rangle}

\renewcommand{\b}{\beta}
\newcommand{\ind}{\mathbf 1}

\numberwithin{equation}{section}
\numberwithin{theorem}{section}

\begin{document}


\title[On Cameron-Martin Theorem and almost sure global existence]
{On Cameron-Martin Theorem and almost sure global existence}

\author{Tadahiro Oh  and Jeremy Quastel}
%
%
%

\address{
Tadahiro Oh\\
School of Mathematics\\
The University of Edinburgh,
and The Maxwell Institute for the Mathematical Sciences\\
James Clerk Maxwell Building\\
The King's Buildings, Peter Guthrie Tait Road
Edinburgh, EH9 3FD, United Kingdom}

\email{hiro.oh@ed.ac.uk}

\address{
Jeremy Quastel\\
Departments of Mathematics and Statistics\\
University of Toronto\\
40 St. George St\\
Toronto, ON M5S 2E4, Canada,
and
School of Mathematics\\
Institute for Advanced Study\\
Einstein Drive, Princeton\\ NJ 08540\\ USA}
\email{quastel@math.toronto.edu}

\thanks{J.Q.~was supported in part by
Natural Sciences and Engineering
 Research Council of Canada, Canada Council for the Arts, through a  Killam Fellowship, and a grant to the Institute for Advanced Study by the Fund For Math.}



\begin{abstract}
In this note, we
discuss various aspects of invariant measures for nonlinear Hamiltonian PDEs.
In particular, we show
almost sure global existence
for some  Hamiltonian PDEs
with  initial data of the form:~``a smooth deterministic function
$+$
 a rough random perturbation'',
as a corollary to
Cameron-Martin Theorem
and known almost sure global existence results
with respect to Gaussian measures on  spaces of functions.

\end{abstract}

\maketitle

\baselineskip = 15pt

\section{Main results}
\label{SEC:1}

\subsection{Introduction}
\label{SUBSEC:1}

In this note, we discuss almost sure global existence results
for some nonlinear Hamiltonian partial differential equations (PDEs)
as corollaries to Cameron-Martin Theorem \cite{CM}.
In particular, we show almost sure global existence
with initial data of the form
\begin{equation}
 u_0(x;\omega) =  v_0 (x) + \phi(x; \omega),
\label{IV1}
 \end{equation}

\noi
where $v_0$ is a deterministic smooth function
and $\phi(\omega)$ is a random function of low regularity.
On $\T$,
the function $\phi$ is given as\footnote{We drop the factor of $2\pi$ throughout the paper,
when it plays no important role.}
\begin{equation}
\phi(x; \omega) = \sum_{n \in \Z \setminus \{0\}}
\frac{g_n(\omega)}{|n|^\al} e^{inx}
\label{FWa}
\end{equation}

\noi
or
\begin{equation}
\phi(x; \omega) = \sum_{n \in \Z }
\frac{g_n(\omega)}{\jb{|n|^\al}} e^{inx}, \qquad \jb{\,\cdot\,} = (1+|\cdot|^2)^\frac{1}{2},
\label{FWb}
\end{equation}

\noi
where
$\{g_n\}_{n \in \Z}$
is a sequence of independent standard complex-valued Gaussian random variables
on a probability space $(\Omega, \mathcal{F}, P)$.
For both \eqref{FWa} and \eqref{FWb},
we easily see that $\phi$ lies almost surely in $H^{\al-\frac{1}{2}-\eps}(\T)$
for any $\eps > 0$
but not in $H^{\al-\frac{1}{2}}(\T)$.
Note that \eqref{FWa} with $\al = 0$ corresponds to the mean-zero
Gaussian white noise on $\T$:
\begin{equation}
\phi(x; \omega) = \sum_{n \in \Z \setminus \{0\}}
g_n(\omega) e^{inx}.
\label{FWc}
\end{equation}

Let us describe one of the motivations for studying the Cauchy problems
with initial data of the form \eqref{IV1},
namely
\begin{equation}
\text{ a smooth deterministic  function }+ \text{ a  rough random perturbation}.
\label{IV2}
\end{equation}

\noi
Given
smooth physical data in an ideal situation,
we may introduce rough and random perturbations to these data
due to the limitations of accuracy in physical observations and storage of such data.
Hence, we believe that it is important
to study  Cauchy problems with initial data of the form \eqref{IV2}.
Initial data \eqref{IV1} with \eqref{FWa} or \eqref{FWb}
are the simplest models for \eqref{IV2}
with rough Gaussian perturbations.
One typical random
noise we introduce in this kind of situation is the white noise,
which appears ubiquitously in the physics literature.
The white noise, however, is very rough
and we can
 handle a smooth initial condition perturbed by the white noise
only in a limited case.

\subsection{Invariant Gibbs measures for Hamiltonian PDEs}
\label{SUBSEC:1.1}
Given a Hamiltonian flow on $\R^{2n}$:
\begin{equation} \label{HR2}
\begin{cases}
\dot{p}_j = \frac{\partial H}{\partial q_j} \\
\dot{q}_j = - \frac{\partial H}{\partial p_j}
\end{cases}
\end{equation}

\noi
 with Hamiltonian $ H (p, q)= H(p_1, \cdots, p_n, q_1, \cdots, q_n)$,
Liouville's theorem states that the Lebesgue measure
$\prod_{j = 1}^n dp_j dq_j$
on $\mathbb{R}^{2n}$ is invariant under the flow.
Then, it follows from the conservation of the Hamiltonian $H$
that  the Gibbs measures
$e^{-\beta H(p, q)} \prod_{j = 1}^{n} dp_j dq_j$ are invariant under the dynamics of \eqref{HR2},
where $\beta> 0$ is the reciprocal temperature.

In the context of the
nonlinear Schr\"odinger equations (NLS) on $\T$:
\begin{equation}\label{NLS}
i u_t - u_{xx} \pm |u|^{p-2}u = 0, \quad (x, t) \in \T\times \R
\end{equation}

\noi
with the Hamiltonian:
\begin{equation}\label{Hamil1}
H(u) = \frac{1}{2}\int_\T |u_x|^2 dx \pm \frac{1}{p}\int_\T |u|^p dx,
\end{equation}

\noi
Lebowitz-Rose-Speer \cite{LRS}
considered the Gibbs measure of the form\footnote{Throughout the paper,
$Z$ denotes various normalizing constants.}:
\begin{align} \label{Gibbs}
d \mu & = d\mu^\beta = Z^{-1} e^{-\beta H(\phi)}  d \phi \notag \\
 & = Z^{-1}
 e^{ \mp \frac{\beta}{p}\int_\T |\phi|^p dx}
 e^{-\frac{\beta}{2}\int_\T |\phi_x|^2 dx} d \phi.
\end{align}

\noi
Here, $d\phi$ denotes the non-existent Lebesgue measure
on the infinite dimensional phase space of functions on $\T$,
and thus the expression \eqref{Gibbs} is merely formal at this point.

Noting that
$ e^{-\frac{\beta}{2}\int |\phi_x|^2 dx} d \phi$
is the Wiener measure on $\T$
with variance $\beta^{-1}$,
Lebowitz-Rose-Speer
showed
 that such a Gibbs measure $\mu$ is
a well-defined probability measure on $H^{\frac{1}{2}-\eps}(\T)$, $\eps > 0$.
In the focusing case, i.e.~with the minus sign in \eqref{Hamil1},
this construction  holds only for $p < 6$
with the $L^2$-cutoff $\ind_{\{\|\phi\|_{L^2} \leq B \}}$ for any $B>0$,
and for $ p = 6$ with sufficiently small $B > 0$.

 Bourgain \cite{BO4} continued this study and
proved invariance of the Gibbs measure $\mu$ under the dynamics of NLS
\eqref{NLS}.\footnote{In order to avoid the problem
at the zero frequency,
we need to insert $-\frac{\beta}{2} \int |u|^2 dx$ in \eqref{Gibbs} for NLS.
As this is standard, we omit this term in the following for simplicity of the presentation.}
See also McKean \cite{MK1} for the cubic case.
The main difficulty in \cite{BO4} was to establish
the global dynamics
 almost surely on the statistical ensemble.
Bourgain achieved this goal
by exploiting invariance
of the finite dimensional Gibbs measures
for the finite dimensional approximations
to \eqref{NLS}.
In the same paper,
he also considered
the generalized KdV equations (gKdV):
\begin{equation}\label{KdV}
 u_t + u_{xxx} \mp u^{p-2}u_x = 0, \quad (x, t) \in \T\times \R.
\end{equation}

\noi
with the Hamiltonian:
\begin{equation}\label{Hamil2}
H(u) = \frac{1}{2}\int_\T u_x^2 dx \pm \frac{1}{p(p-1)}\int_\T u^p dx.
\end{equation}

\noi
In particular,
invariance of the Gibbs measures
for KdV ($p = 3$) and mKdV ($p = 4$)
was established in \cite{BO4}.
Recently,
Richards \cite{Richards} treated the case of the quartic KdV ($p = 5$).
There have been
papers in this direction
by Bourgain
\cite{BO5,  BO6, BO96, BO8, BO9}
and and other mathematicians that followed his idea
 \cite{TZ1, TZ2, BT1, BT3,
OH3, OHSBO, TZ3, NORS, dS1, dS2, Deng}.

In the following, we set $\beta = 1$ for simplicity.
Then, the Gibbs measure $\mu$ in \eqref{Gibbs}
is absolutely continuous with respect to
 the Wiener  measure $\rho$ with the density:
\begin{equation} \label{Gauss0}
 d\rho = Z^{-1} e^{-\frac{1}{2}\int_\T |\phi_x|^2 dx} d \phi.
 \end{equation}

\noi
A typical element $\phi$ in the support of the Wiener measure
can be represented by the Fourier-Wiener series \eqref{FWa} with $\al = 1$.

In the defocusing case, the Gibbs measure $\mu$
and the Wiener measure $\rho$ are equivalent,
 i.e.~mutually absolutely continuous.
In particular,
almost sure global existence with respect to the Gibbs measure $\mu$
implies
almost sure global existence with respect to the Wiener measure $\rho$.
For example, the defocusing NLS \eqref{NLS} for any $p$
is almost surely globally well-posed
with respect to the random initial data
$u|_{t = 0} = \phi$, where $\phi$ is as in \eqref{FWa} with $\al = 1$.
We point out that, for $p> 6$,
 this is beyond the known deterministic global well-posedness
results.
In Subsection \ref{SUBSEC:3},
we show how this result
can be extended to almost sure global
well-posedness
for the initial data $v_0 +\phi$,
where $v_0 \in H^1(\T)$ and $\phi $ is as in \eqref{FWa} with $\al = 1$.

\begin{remark}\rm

In the following, we recall two properties of Gibbs measures.
Although they are well known in probability theory and in statistical mechanics,
we decided to include this remark for readers' convenience,
in particular, for those in PDEs.

\smallskip

\noi
(i) Variational characterization of the Gibbs measure.
Here, we restrict our attention to the finite dimensional setting \eqref{HR2}.
With $\phi = (p, q)$,
 the Gibbs measure
can be written as
\begin{equation}
 d \mu_\beta  = f_\beta^*(\phi) d\phi : = Z_\beta ^{-1} e^{-\beta H(\phi)} d\phi ,
\label{X1}
 \end{equation}

\noi
where $d\phi$ denotes the Lebesgue measure $d\phi  = \prod_{j = 1}^n dp_j dq_j$ on $\R^{2n}$.

Given a probability measure $\rho$
that is absolutely continuous with respect to the Lebesgue measure $d\phi$,
we define its entropy $S(\rho)$ and average energy $\jb H(\rho)$ by
\[S(\rho) = -\int \frac{d \rho}{d \phi} (\phi) \log \Big(\frac{d \rho}{d \phi}(\phi)\Big) d\phi
\qquad \text{and}
\qquad \jb{ H}(\rho) = \int H(\phi) \frac{d \rho}{d \phi}(\phi) d\phi,\]

\noi
respectively, where $H$ is the Hamiltonian for the underlying dynamics.
In the following, we consider the maximization
problem of the entropy $S(\rho) $ for a given average energy $\jb H(\rho) = C$.
We  assume that, for a given value of $C$,
there exists a unique $\beta > 0$ such that
$\jb H(\mu_\beta) = C$.
For simplicity of notations,
we write $S(f)$ and $\jb{H}(f)$ for $S(\rho)$ and $\jb H(\rho)$,
where $f : = f_\rho = \frac{d\rho}{d \phi}$
denotes the Radon-Nikodym derivative of $\rho$ with respect to the Lebesgue measure $d \phi$.
Then, by the Lagrange multiplier method
with two constraints $\jb H(f) = C$
and $M(f) := \int f (\phi) d\phi = 1$,
we have
\begin{align*}
 dS(f) & = \be d \jb H(f) + \g d M(f)\\
& \LRA  \int   \big( \log f(\phi) + 1 + \g + \be H(\phi)\big) g(\phi) d\phi = 0
\end{align*}

\noi
for all test functions $g$.
Thus, we conclude that
$f(\phi) = e^{-1 -\g -\be H(\phi)}$.
Moreover, by the mass constraint $M(f) = 1$,
we must have
$f(\phi) = Z_\be^{-1} e^{ -\be H(\phi)} = f^*_\be(\phi)$,
where $f^*_\be$ is as in \eqref{X1}.
Hence, if there is any extremal point
for the entropy functional,
it has to be the  Gibbs measure $\mu_\be$.
Also, by a direct computation, we have
$d^2 S(f)(g, g) = - \int \frac{g^2}{f} d\phi \leq 0$.
Therefore,  the Gibbs measure
 $\mu_\be$ is the unique maximizer of
 the entropy for a given average energy.

\medskip

\noi
(ii) Dependence of the Gibbs measure $\mu$ on $\beta >0$:
In the mathematics literature, the value of $\beta$
is often set to be $1$ for simplicity.
In the following, we
discuss the relation of  $\mu^\beta$ for different values of  $\beta > 0$.
In particular, we show that the Gibbs measures
$\mu^\beta$ and $\mu^\g$, $\beta, \g > 0$,
are singular if $\beta \ne \g$.

Consider the Gaussian measure $\rho^\beta$ with the density:
\[ d\rho^\beta = Z^{-1} e^{-\frac{\beta}{2}\int_\T |\phi_x|^2 dx} d \phi.\]

\noi
This is a Gaussian probability measure
on $\dot H^s(\T)$, $s< \frac{1}{2}$.
Indeed, with $B_\beta := \beta^{-1} D^{2s-2}$,
where $D = \sqrt{-\Dl}$,
we have
\begin{align*}
-\frac{\beta}{2}\int_\T|\phi_x|^2 dx
= -\frac{1}{2}\jb{B_\beta^{-1} \phi, \phi}_{\dot H^s}.
\end{align*}
	
\noi
Hence, $\rho_\beta$ is the (mean-zero) Gaussian
measure with the covariance operator $B_\beta$.
Moreover, with $e_n:= |n|^{-s} e^{inx}$, $n \ne 0$,
we have
$ B_\beta e_n = \ld_n(\beta) e_n$,
where
\begin{equation}
\ld_n (\beta)= \beta^{-1} |n|^{2s-2}.
\label{X2}
\end{equation}

\noi
Now, consider two Gaussian measures
$\rho^\beta$ and $\rho^\g$,
   $\beta, \g > 0$.
Feldman-H\'ajek theorem \cite{Feldman, Hajek}
states that two Gaussian measures
are either (i) equivalent 
or  (ii) singular.
Moreover, 
letting
\[ S(\beta, \g) = \sum_{n \ne 0} \frac{(\ld_n(\beta) - \ld_n(\g))^2}
{(\ld_n(\beta) + \ld_n(\g))^2}, \]

\noi
we have
(i) $\rho^\beta$ and $\rho^\g$
are equivalent if $S(\beta, \g) < \infty$
and
(ii) $\rho^\beta$ and $\rho^\g$
are singular if $S(\beta, \g) = \infty$.
See also Kakutani's dichotomy theorem \cite{Kakutani}
on equivalence 
of infinite product measures.

With \eqref{X2}, it is easy to see that $S(\b, \g) = \infty$
if $\b \ne \g$.
Thus, $\rho^\b$ and $\rho^\g$ are singular, if $\beta \ne \g$.
Therefore,
noting that $\mu^\beta$ is absolutely continuous
with respect to $\rho^\beta$,
it follows that the Gibbs measures
$\mu^\beta$ and $\mu^\g$ are singular if $\b \ne \g$.

\end{remark}

Next, recall that NLS \eqref{NLS}
and gKdV \eqref{KdV} also preserve
the $L^2$-norm of solutions.
Thus, the Gaussian white noise $\mu_0$ with the density:
\begin{equation} d\mu_0 = Z^{-1} e^{-\frac{1}{2}\int_\T |\phi|^2 dx} d \phi
\label{white}
\end{equation}

\noi
is expected to be invariant for these equations.
On $\T$, a typical element $\phi$ in the support of the white noise $\mu_0$
is represented
by \eqref{FWc} (in the mean-zero case),
which is almost surely in $H^{-\frac{1}{2}-\eps}(\T)$ for any $\eps > 0$
but not in $H^{-\frac{1}{2}}(\T)$.
It is this low regularity that makes it difficult to
rigorously study invariance of the white noise.
Nonetheless,
for KdV \eqref{KdV} with $p = 3$,
the (mean-zero) white noise is shown to be
invariant \cite{QV, OH4, OHRIMS, OQV}.
See also \cite{OH6, OH7}.
In particular, this result yields
almost sure global existence for KdV with the white noise as initial data,
namely with $u|_{t = 0} = \phi$, where $\phi$ is as in \eqref{FWc}
conditioned that $g_{-n} = \cj{g_n}$.

Note that
almost sure existence of  a solution with the white noise as initial data
(but not its invariance) also follows from
deterministic global well-posedness of KdV in $H^{-1}(\T)$
by Kappeler-Topalov \cite{KT}, exploiting the integrable structure of the equation.
However, the result in \cite{OH4, OHRIMS}
can be applied to non-integrable variants of KdV,
and moreover it asserts a stronger form of uniqueness.

In \cite{OQV},
the white noise was shown to be a weak limit
of invariant measures, more precisely, a limit
of interpolations of the Gibbs measures (with a parameter) and the white noise.
This result holds not only for KdV but also for  cubic NLS and mKdV,
i.e.~\eqref{NLS} and \eqref{KdV} with $p = 4$.
Due to lack of well-defined dynamics in the support
of the white noise,
this does not yield invariance of the white noise
for cubic NLS and mKdV,
but it only
 provides  a strong evidence of such invariance.

\subsection{Probabilistic Cauchy theory}
\label{SUBSEC:2}

In an effort to study the Cauchy problem
for cubic NLS in low regularity,
Colliander-Oh \cite{CO}
considered the following
 {\it Wick ordered cubic NLS} on $\T$:
\begin{equation}
	\label{NLS1}
\textstyle
		i u_t - u_{xx} \pm u (|u|^2 -2  \fint |u|^2 dx) = 0,
\end{equation}

\noi
with random initial data of the form \eqref{FWb},
where
$\fint |u|^2 dx := \frac{1}{2\pi} \int |u|^2 dx$.
This equation first appeared in \cite{BO96}
in the context of the defocusing cubic NLS on $\T^2$
as an equivalent formulation of
the Hamiltonian equation arising from
 the Wick ordered Hamiltonian.

Note that
 $u$ solves \eqref{NLS} if and only if
$v(t) = e^{i \g t} u(t)$, with $\g \in \R$,
 solves $i
\partial_t v - v_{xx} \pm |v|^2v + \gamma v = 0$.
Hence, by letting $\g= \mp 2 \fint |u|^2 dx$ along with the $L^2$-conservation,
\eqref{NLS} is equivalent to \eqref{NLS1},
at least  for $u_0 \in L^2(\T)$.
For $u_0 \notin L^2(\T)$,
we cannot freely convert solutions of \eqref{NLS1} into solutions of \eqref{NLS}.
See \cite{OS} for more discussions on the relation between the cubic NLS \eqref{NLS} and
the Wick ordered cubic NLS \eqref{NLS1}.

In \cite{CO},
it is shown that
\eqref{NLS1} is
almost surely
locally  well-posed with the initial data \eqref{FWb} with $\al > \frac{1}{6}$,
corresponding to $H^s(\T)$, $s > -\frac{1}{3}$,
and
almost surely
globally  well-posed with the initial data \eqref{FWb} with $\al > \frac{5}{12}$,
corresponding to $H^s(\T)$, $s > -\frac{1}{12}$.
Note that $\phi$ in \eqref{FWb}
represents a typical element
of  the following Gaussian measure $\rho_\al$ with the density:
\begin{equation}
d \rho_\al = Z e^{-\frac{1}{2}\int_\T |\phi|^2 dx - \frac{1}{2}\int_\T |D^\al \phi|^2 dx }d\phi.
\label{Gauss1}
\end{equation}

The probabilistic local argument in \cite{CO} closely follows
that by Bourgain \cite{BO96}.
The main ingredients are (i) an improvement of the Strichartz estimates
under randomization of  initial data
and (ii) hypercontractivity of the Ornstein-Uhlenbeck semigroup.
Burq-Tzvetkov \cite{BT2} exploited (i)
to establish an almost sure local existence result for the nonlinear wave equation (NLW)
for a wider class of randomizations.

The probabilistic global argument in \cite{CO}
was the first almost sure global existence argument in the absence of conservation laws
or formally invariant measures.
The proof was based on the adaptation of Bourgain's high-low method \cite{BO98}
in the probabilistic setting.
In particular, it exploited
the $L^2$-conservation
and invariance of the Gaussian measure $\rho_\al$ in \eqref{Gauss1}
under the {\it linear} flow.
More recently,
there have been almost sure global existence results
for some other equations
\cite{BT4, NPS, BTT} in the absence of
conservation laws
or formally invariant measures.
The argument is based on a combination
of a conservation law at a higher regularity
and a probabilistic argument.

\medskip

In the following, we focus on the almost sure global existence result in \cite{CO}.
It says that
 given $\al > \frac{5}{12}$,
 there exists $\Sigma_0 = \Si_0(\al)$ with $\rho_\al(\Si_0) = 1$
 such that, if $\phi \in \Si_0$,
 then there exists a global solution $u$ to \eqref{NLS1} with $u|_{t = 0} = \phi$.
 There are two issues about this almost sure global existence result:
 \begin{itemize}
\item
It does not say anything about what happens to $\Si_0$ under the dynamics.
In particular, it does {\it not} guarantee that $\Si_0$ remains a set of full measure under the \eqref{NLS1} flow.
Let $\Phi(t) : \phi \mapsto u(t) = \Phi(t) \phi$
be the solution map of \eqref{NLS1}.
Then, it may happen that
$\Phi(t) \Si_0$ for $t>0$ is a set of smaller measure
and we may even have $\rho_\al \big(\Phi(t) \Si_0\big) = 0$
for some $t>0.$

\item
The uniqueness statement for the local result in \cite{CO}
states the following;
if $\phi = \phi(\omega)$ is a ``good'' initial condition,
then the solution $u(t) = \Phi(t) \phi$
exists up to time $\dl > 0$
and  uniqueness holds in the ball
centered at $S(t) \phi$ of radius $1$ in $X^{0, \frac{1}{2}+}_{[0, \dl]}$.
Here,
$S(t) = e^{-i t \dx^2}$ denotes the linear propagator for \eqref{NLS1}
and
$X^{0, \frac{1}{2}+}_{[0, \dl]}$
denotes the local-in-time
version of the $X^{s, b}$ space
onto the time interval $[0, \dl]$ (with $s = 0$ and $b = \frac{1}{2}+$).
This is a typical uniqueness statement for
the probabilistic local Cauchy theory.  See \cite{BO96, BT2}.
However, the uniqueness statement for the almost sure global existence
result in \cite{CO} holds in a much milder sense.
See Remark 1.2  
in \cite{CO}.
 \end{itemize}

\noi
The next theorem addresses both of the issues described above.
\begin{theorem} \label{THM:NLS}
Let $\al > \frac{5}{12}$.
Then, there exists a set $\Si \subset H^{\al - \frac{1}{2}-\eps}(\T)$, $\eps > 0$, of full measure with respect to $\rho_\al$
such that
\begin{itemize}
\item[\textup{(i)}] $\Si$ is invariant under the \eqref{NLS1}-dynamics.
In particular, $\rho_\al(\Phi(t)\Sigma) = 1$ for any $t \in \R$
and if $\phi \in \Si$, then the corresponding solution
$u(t) = \Phi(t) \phi$ exists globally.

\item[\textup{(ii)}]
Given $\phi \in \Si$,
the global solution $u(t) = \Phi(t) \phi$  is unique in the following sense.
Given $t_* \in \R$,
there exists positive $\dl = \dl (\phi, t_*)>0$
such that uniqueness holds in the ball centered at $S(\cdot - t_*)u(t_*)$
of radius 1 in $X^{0, \frac{1}{2}+}_{[t_* - \dl, t_*+ \dl]}$.
Moreover,
for each finite time interval $I$,
 $\dl>0$ is bounded away from 0
 for all $t_* \in I$.
\end{itemize}
\end{theorem}

\noi
We point out that this uniqueness statement
is in the spirit of the usual probabilistic local Cauchy theory
and is stronger than the uniqueness statement for almost sure global solutions in \cite{CO}.
We present
the proof of Theorem
\ref{THM:NLS}
in  Section \ref{SEC:NLS}.
Previously, Burq-Tzvetkov \cite{BT4}
 constructed an invariant set
of full measure
in considering almost sure global existence for NLW on $\T^3$.
Their idea was based on first characterizing the set
of initial data
such that
the corresponding linear solutions
satisfy some space-time bounds,
guaranteeing global existence,
and then showing that random initial data almost surely belongs to this set.
The global argument in \cite{CO} exploits
 finer properties of {\it products} of the linear solutions with random initial data
 (such as the hypercontractivity of the Ornstein-Uhlenbeck semigroup),
 which is difficult to characterize in terms of  individual initial data.
Hence,  the proof of Theorem \ref{THM:NLS} follows
 a different path than that in \cite{BT4}.

\subsection{Cameron-Martin Theorem
and almost sure global existence}
\label{SUBSEC:3}

In this section, we recall Cameron-Martin Theorem
and discuss its implications in the context of almost sure global existence
for nonlinear Hamiltonian PDEs.

For this purpose, we first need to briefly go over the definition
of abstract Wiener spaces introduced by
Gross \cite{GROSS}.
See also
Kuo \cite{KUO}.
Let $H$ be a real separable Hilbert space.
It is known that the Gauss measure
$\rho$ with the density
$d \rho = Z^{-1} e^{-\frac{1}{2}\|x\|_{H}^2}dx$
is only finitely additive if $\dim H = \infty$.

Let $\mathcal{P}$ denotes the collection
of all finite dimensional orthogonal projections of $H$.
A seminorm $|||\cdot|||$
on $H$ is said to be measurable
if, for any $\eps > 0$, there exists $\mathbb{P}_\eps \in \mathcal{P}$
such that
$\rho (|||\mathbb{P}x|||>\eps) < \eps$
for all $\mathbb{P}\in \mathcal{P}$ with $\mathbb{P}\perp \mathbb{P}_\eps$.
Let $B$ be the completion of $H$ with respect to this seminorm $|||\cdot|||$.
Then,
Gross \cite{GROSS} showed that
$\rho$ can be made sense of as  a countably additive
Gaussian measure on $B$.
In this case, we say that
the triplet $(B, H, \rho)$ is an {\it abstract Wiener space}.
The original Hilbert space $H$
is often referred to as a {\it Cameron-Martin space}
or a reproducing kernel Hilbert space.

Let $(B, H, \rho)$ be an abstract Wiener space.
Then,  Cameron-Martin Theorem
states the following.

\smallskip

\noi
{\bf Cameron-Martin Theorem.}
{\it Given $h \in B$, define
a shifted measure $\rho_h$ by $ \rho_h(\, \cdot\,) := \rho(\, \cdot\,  - h)$.
Then,
the shifted measure $\rho_h$
is mutually absolutely continuous with respect to $\rho$
if and only if
$h \in H$.}

\smallskip

\noi
This theorem also provides
a precise expression of the Radon-Nikodym derivative.
This absolute continuity under a shift in the direction of $H$
leads to the $H$-differentiation,
which
plays a key role in the Malliavin Calculus.
See Shigekawa \cite{Shigekawa}.

\begin{example}\rm
Consider the Wiener measure $\rho$ in \eqref{Gauss0}.
More precisely, consider the Gaussian measure $\rho$
with the density:
\begin{equation}
  d\rho = Z^{-1} e^{-\frac{1}{2}\int_\T |\phi|^2 dx-\frac{1}{2}\int_\T |\phi_x|^2 dx} d \phi
= Z^{-1} e^{-\frac{1}{2}\|\phi\|_{H^1}^2} d\phi.
\label{Gauss2}
\end{equation}

\noi
Then, $\rho$
is the Gauss measure on $H = H^1(\T)$.
It is known that,  with $B = H^s(\T)$, $s < \frac{1}{2}$,
the triplet $(B, H, \rho)$ is an abstract Wiener space.
See B\'enyi-Oh \cite{Benyi}
for examples of other Banach spaces $B$
such that $(B, H, \rho)$ is an abstract Wiener space.
Note that (i) this Gaussian measure $\rho$ is absolutely continuous with respect to
the Gibbs measure $\mu$ in the defocusing case, i.e.~with the minus sign in \eqref{Gibbs}
and (ii) the Fourier-Wiener series \eqref{FWb} represents functions in the support of $\rho$.
Then, as a corollary to
invariance of the Gibbs measure $\mu$ and Cameron-Martin Theorem,
we have the following statement.

\begin{theorem} \label{THM:1}
Let $v_0 \in H^1(\T)$.
Then, the solution $u = u( x, t; \omega)$ to the defocusing NLS \eqref{NLS}
with the initial data of the form
\[ u_0(x;\omega) = v_0(x) + \sum_{n\in \Z} \frac{g_n(\omega)}{\jb{n}} e^{inx}
\in H^{\frac{1}{2}-}(\T) \setminus H^{\frac{1}{2}}(\T), \text{ a.s.}, \]

\noi
exists globally in time, almost surely in $\omega$.
\end{theorem}

\noi
For the cubic and quintic NLS ($p = 4$ and $p = 6$, respectively),
Theorem \ref{THM:1} follows from the deterministic global well-posedness
results by Bourgain \cite{BO1, BO04}.
For $p > 6$,
Theorem \ref{THM:1} does not follow from known (deterministic) results.

In the focusing case,
the Gibbs measure $\mu$ comes with an $L^2$-cutoff
$\ind_{\{\|\phi\|_{L^2} \leq B\}}$.
Since a shift by $v_0$ does not preserve
this $L^2$-cutoff, a result analogous to Theorem \ref{THM:1}
does not hold as a corollary to Cameron-Martin Theorem.
Recall, however,  that
the Gibbs measure makes sense only for $p\leq 6$ in the focusing case
(where $B$ is sufficiently small if $p = 6$),
where deterministic global well-posedness results are available in
the regularity of the Gibbs measures.
When $p = 4$, the cubic NLS \eqref{NLS} is globally well-posed in $L^2(\T)$.
When $ p = 6$,
a modification of the argument in  \cite{BO04}
yield global well-posedness of the quintic NLS for data with small $L^2$-norms.

\end{example}

\begin{example}\rm
In this example, we assume that all the functions are real-valued
with mean zero on $\T$.
Let $\mu_0$ be the mean-zero Gaussian white noise defined in \eqref{white}.
Note that $\mu_0$ is the Gauss measure on $H = L^2_0(\T)$,
where $L^2_0(\T)$ denotes the collection of real-valued
functions in $L^2(\T)$ with mean zero on $\T$.
With $B = H^s(\T)$, $s < -\frac{1}{2}$,
the triplet $(B, H, \mu_0)$ forms an abstract Wiener space.
Hence,
as a corollary to invariance of the white noise for KdV \eqref{KdV} with $p = 3$
and Cameron-Martin Theorem, we have the following.

\begin{theorem} \label{THM:2}
Let $v_0 \in L^2_0(\T)$.
Then, the solution $u = u( x, t; \omega)$ to KdV \eqref{KdV} with $p = 3$
with the initial data of the form
\[ u_0(x; \omega) = v_0(x) + \sum_{n\in \Z \setminus \{0\}} g_n(\omega) e^{inx}
\in H_0^{-\frac{1}{2}-}(\T) \setminus H_0^{-\frac{1}{2}}(\T), \text{ a.s.}, \]

\noi
exists globally in time, almost surely in $\omega$.
Here,
$H^s_0(\T)$
 denotes the collection of real-valued
functions in $H^s(\T)$ with mean zero on $\T$.
\end{theorem}

\noi
Note that Theorem \ref{THM:2}
also follows from the deterministic global well-posedness results
by Kappeler-Topalov \cite{KT}.
However,  the result in \cite{KT}
is not applicable to non-integrable variants of KdV,
while Theorem \ref{THM:2}
also holds for some non-integrable variants of KdV.


\end{example}

\begin{example}\rm
The Gaussian measure $\rho_\al$ defined in \eqref{Gauss1}
 is the Gauss measure on $H = H^\al(\T)$.
With $B = H^s(\T)$, $s < \al- \frac{1}{2}$,
the triplet $(B, H, \rho_\al)$  forms an abstract Wiener space.
Hence,
as a corollary to Theorem 2 in \cite{CO}
and Cameron-Martin Theorem, we have the following.

\begin{theorem} \label{THM:3}
Let $v_0 \in H^\al(\T)$, $\al > \frac{5}{12}$.
Then, the solution $u = u( x, t; \omega)$ to the Wick ordered cubic NLS \eqref{NLS1}
with the initial data of the form
\begin{equation}
 u_0(x; \omega) = v_0(x) + \sum_{n\in \Z } \frac{g_n(\omega)}{\jb{|n|^\al}} e^{inx}
\in H^{\al -\frac{1}{2}-}(\T) \setminus H^{\al -\frac{1}{2}}(\T), \text{ a.s.},
\label{IV3}
 \end{equation}

\noi
exists globally in time, almost surely in $\omega$.

\end{theorem}

\noi
It follows from a slight modification of the proof of Theorem 1 in \cite{CO}
 that the solution $u$ to \eqref{NLS1}
with initial data \eqref{IV3}, $\al > \frac{1}{3}$,
exists locally in time, almost surely in $\omega$,
 even if $v_0$ is only in $L^2(\T)$.
However, it seems that  much more effort is required
to modify the global argument
in \cite{CO} to obtain Theorem \ref{THM:3}
for $v_0 \in L^2(\T)$.

\end{example}

\subsection{On absolute continuity under a shift for other classes of randomizations}
There are several results on
almost sure global existence
for a more general class of randomized initial data.
See \cite{BT4, NPS, BTT}.
In this subsection, we discuss
the effect of a shift by a smooth function on  such randomized initial data.
For simplicity, we restrict our attention to $\T^d$.
Fix a function $u = \sum_{n \in \Z^d} \ft u_n e^{in\cdot x} $ in $H^s(\T^d)$,
and define its randomization 
$u^\omega$ by
\[ u^\omega = \sum_{n \in \Z^d} a_n(\omega) \ft u_n e^{in\cdot x},\]

\noi
where $\{a_n\}_{n \in \Z^d}$
is a sequence of independent 
complex-valued random variables
on a probability space $(\Omega, \mathcal{F}, P)$.
Given a deterministic function $v$ on $\T^d$,
we consider a shifted function
\[ w^\omega = v + u^\omega
= \sum_{n \in \Z^d} \big(\ft v_n + a_n(\omega) \ft u_n\big) e^{in\cdot x}.\]

\noi
We would like to know when
the distribution
of $w^\omega$ is absolutely continuous with respect to
that of $u^\omega$.
Clearly, the support of $a_n(\omega) \ft u_n$
must contain the support of $\ft v_n + a_n(\omega) \ft u_n$.
This eliminates a certain class of random variables
such as the Bernoulli random variables.
Moreover,
since our interest is to
determine a class of functions $v$
such that
the distribution
of the shifted random function $w^\omega$ is absolutely continuous with respect to
that of $u^\omega$,
we may assume that
$\{a_n(\omega)\ft u_n\}_{\omega \in \Omega} = \C$ for each $n \in \Z^d$.
For simplicity,
we set $\{a_n\}_{n\in \Z^d}$ to be
a sequence of independent standard
complex-valued Gaussian random variables
$\{g_n\}_{n\in \Z^d}$ in the following.

Recall the following definition of the Hellinger integral.
See, for example,  Da Prato \cite{DaPrato}.
Given two probability measures $\mu$ and $ \nu$,
the Hellinger integral of $\mu$ and $\nu$ is defined by
\[ H(\mu, \nu) = \int_\Omega \sqrt{\frac{d \mu}{d \zeta}
\frac{d \nu}{d \zeta}} d\zeta,\]

\noi
where $\zeta = \frac{1}{2}(\mu + \nu)$.
Note that $\mu$ and $\nu$ are absolutely continuous
with respect to $\zeta$, so the Radon-Nikodym derivatives
$\frac{d \mu}{d \zeta}$ and $\frac{d \nu}{d \zeta}$ make sense.
If $\nu$ is absolutely continuous
with respect to $\mu$,
we can write the Hellinger integral as
\[ H(\mu, \nu) = \int_\Omega \sqrt{\frac{d \nu}{d \mu}} d\mu.\]

\noi
Given $\{\mu_n\}_{n \in \mathbb{N}}$ and $\{\nu_n\}_{n \in \mathbb{N}}$ are sequences of
probability measures
on $\C$, consider the product measures on $\C^\infty$:
$ \mu = \bigotimes_{n = 1}^\infty \mu_n$
and $ \nu = \bigotimes_{n = 1}^\infty \nu_n$.
In this case, the Hellinger integral of $\mu$ and $\nu$
is given by
$ H(\mu, \nu) = \prod_{n = 1}^\infty H(\mu_n, \nu_n)$.
Then, Kakutani's theorem \cite{Kakutani} states that
(i) $\mu$ and $\nu$ are equivalent
if $H(\mu, \nu) > 0$,
and (ii)
$\mu$ and $\nu$ are singular
if $H(\mu, \nu) = 0$.

Now, let $\mu_n$ and $\nu_n$
be the  probability measures on $\C$
induced by
$\omega \mapsto g_n(\omega) \ft u_n$
and
$\omega \mapsto \ft v_n + g_n(\omega) \ft u_n$,
respectively.
Namely, the density functions are given by
\[ d\mu_n = \frac{1}{2\pi} e^{-\frac{1}{2} \frac{|z|^2}{|\ft u_n|^2} }dz
\quad \text{and}\quad
 d\nu_n = \frac{1}{2\pi} e^{-\frac{1}{2} \frac{|z - \ft v_n |^2}{|\ft u_n|^2} }dz
 =
 e^{-\frac{1}{2} \frac{|\ft v_n |^2 - \text{Re}\jb{\ft v_n, z}  }{|\ft u_n|^2} }
d\mu_n.
 \]
Then, the Hellinger integral of $\mu_n$ and $\nu_n$
is given by
$ H(\mu_n, \nu_n) = e^{-\frac{1}{8}\frac{|\ft v_n|^2}{|\ft u_n|^2}}$.
Let $ \mu = \bigotimes_{n \in \Z^d}^\infty \mu_n$
and $ \nu = \bigotimes_{n \in \Z^d}^\infty \nu_n$.
Then, $\mu$ and $\nu$ represent
the probability distributions of (the Fourier coefficients of) $u^\omega$
and $w^\omega$, respectively.
Moreover, we have
\[ H(\mu, \nu) = \prod_{n \in \Z^d}  e^{-\frac{1}{8}\frac{|\ft v_n|^2}{|\ft u_n|^2}}.\]

\noi
Hence,
$\mu$ and $\nu$ are equivalent
if and only if $H(\mu, \nu)^{-1} < \infty$,
i.e.
\begin{equation}
\sum_{n\in \Z^d}\frac{|\ft v_n|^2}{|\ft u_n|^2 }< \infty.
\label{sum}
\end{equation}

\noi
In particular, if $\ft u_n = 0$ for some $n$, we must have $\ft v_n = 0$.

In general, given $u \in H^s(\T^d)$,
it may not be easy to determine a class of functions $v$ such that \eqref{sum}
holds.
For example, suppose that
$\ft u_n = |n|^{-1}$, $n \in \Z^d \setminus \{0\}$,
and $\ft u_0 = 0$.
Namely, $u^\omega$ is the mean-zero Gaussian free field on $\T^d$.
In this case, we have $u \in \dot H^s(\T^d)\setminus \dot H^{1-\frac{d}{2}}(\T^d)$,
$s < 1-\frac{d}{2}$, almost surely.
Since the randomization
on the Fourier coefficients does not introduce any smoothing
in terms of differentiability almost surely in $\omega$,
we also have $u^\omega
\in \dot H^s(\T^d)\setminus \dot H^{1-\frac{d}{2}}(\T^d)$,
$s < 1-\frac{d}{2}$, almost surely.
In view of the  condition \eqref{sum},
we see that
the distributions of $u^\omega$
and the shifted random variable $w^\omega$ are equivalent
if $v \in \dot H^1(\T^d)$,
and that they are singular
if $v \not\in \dot H^1(\T^d)$.
If one knows what kind of noise $u^\omega$ is added to smooth initial data $v$,
then it is possible to repeat the computation above.

\subsection{On the large deviation principle
with respect to small random perturbations}

In this subsection,
we discuss the large deviation principle
for solutions with
 initial data perturbed by small random noises.
 In particular,
 we consider initial data
 of the form:
\begin{equation}
 u_0^\eps(x;\omega) =  v_0 (x) + \eps\phi(x; \omega)
\label{Z1}
 \end{equation}

\noi
for small $\eps > 0$,
where $v_0$ is a deterministic smooth function
and $\phi(\omega)$ is as in \eqref{FWa}, \eqref{FWb}, or \eqref{FWc}.
The theory of large deviations
was formalized in the seminal paper by Varadhan \cite{V1}
and we follow his definition.
See also Varadhan \cite{V2}.

In the following, we consider KdV,  \eqref{KdV} with $p = 3$,
and use the notations in Theorem \ref{THM:2}.
Fix  $v_0 \in L^2_0(\T)$
and let $\phi(\omega)$ be the mean-zero Gaussian white noise given by  \eqref{FWc}.
For these $v_0$ and $\phi(\omega)$,
let $u^\eps(\omega)$ be the global solution to KdV 
with initial data $u_0^\eps(\omega)$ defined in \eqref{Z1}.
Note that Theorem \ref{THM:2}
guarantees global existence
of such $u^\eps$ almost surely. 
Fix $s <  -\frac{1}{2}$.
Then, the map $\omega \mapsto u^\eps(\omega)$
induces probability measures $\mu_\eps$
on $C(\R; H^{s}_0(\T))$.
In the following, we discuss the large deviation principle for $\mu_\eps$.

First,
we discuss the large deviation principle for
the probability measures
$\rho_\eps = P\circ (u_0^\eps)^{-1}$
on initial data $u_0^\eps$ in \eqref{Z1}.
Define a rate function $I: H^s_0(\T) \to [0, \infty]$
by
\begin{equation}
I(f) = \frac{1}{2} \|f - v_0\|_{L^2_0(\T)}^2.
\label{Z2}
\end{equation}
	
\noi
Note that
(i) $I$ is lower semicontinuous by Fatou's lemma
and (ii) $K_r = \{ f \in H^s_0(\T)  : I(f) \leq r\} \subset L^2_0(\T)$
is compact in $H^s_0(\T)$ for each finite $r \geq 0$.
Then, the large deviation principle
holds for $\{\rho_\eps\}_{\eps > 0}$
with the rate function $I$ in \eqref{Z2}.
Namely, we have
\begin{equation}
\limsup_{\eps \to 0} {\eps^2} \log \rho_\eps (F) \leq - \inf_{f \in F} I(f)
\label{Z3}
\end{equation}

\noi
for any closed set  $F\subset H^s_0(\T)$
and 	
\begin{equation}
\liminf_{\eps \to 0} {\eps^2} \log \rho_\eps (G) \geq - \inf_{f \in G} I(f)
\label{Z4}
\end{equation}

\noi
for any open set $G \subset H^s_0(\T)$.
The inequalities \eqref{Z3} and \eqref{Z4}
follow from Theorems 3.3 and 4.2
in Chapter 3 of Freidlin-Wentzell \cite{FW}.
Note that $K_0 = \{ f \in H^s_0(\T): I(f) = 0\}$
consists of a single function $v_0$.
Hence, it follows
from Remark 2.3 in \cite{V2}
that  $\rho_\eps$ converges weakly to $\dl_{v_0}$
as $\eps \to 0$.

Next, we discuss the large deviation principle
for the probability measures $\mu_\eps$
on solutions to KdV.
First, for fixed $s \in [-1, -\frac 12)$, endow
$C(\R: H^s_0(\T))$
 with  the topology of compact convergence (compact-open topology)
 induced by the usual metric:
 \[ d(u, v) = \sum_{j= 1}^\infty 2^{-j} \frac{\|u-v\|_{L^\infty_t([-j, j]: H^s)}}
 {1+\|u-v\|_{L^\infty_t([-j, j]: H^s)}}.\]

\noi
Then, $C(\R: H^s_0(\T))$
 is a Polish space.
In view of global well-posedness of KdV in $H^{-1}_0(\T)$ by Kappeler-Topalov \cite{KT},
let $X \subset C(\R: H^s_0(\T))$ denote the collection of global-in-time solutions to KdV
constructed in \cite{KT}, endowed with the subspace topology.
It follows from the continuity of the solution map
$\Phi : u(0) \in H^s_0(\T) \mapsto \Phi(u(0)): = u \in X$
with respect to the topology induced by the metric $d(\cdot, \cdot)$
that $X$ is also a Polish space.

Let  $\Psi: = \Phi^{-1}:X \to H^s_0(\T)$ be the evaluation map given
by $\Psi (u) = u(0)$.
By definition of $\mu_\eps$, we have
$\mu_\eps(A) = \rho_\eps(\Psi (A))$.
Now,
define a rate function $\wt I: X \to [0, \infty]$
by
\begin{equation}
\wt I(u) := I(\Psi(u))  = \frac{1}{2} \|\Psi(u) - v_0\|_{L^2_0(\T)}^2.
\label{Z5}
\end{equation}
	
\noi
The lower semicontinuity of $\wt I$
directly follows from that of $I$.
Let $\wt K_r = \{ u \in X  : \wt I(u) \leq r\} \subset C(\R: L^2_0(\T))$.
Given a sequence $\{u_n\}_{n = 1}^\infty \subset \wt K_r$,
it follows from \eqref{Z5} that $\{u_n(0)\}_{n=1}^\infty$
is bounded in $L^2_0(\T)$ and thus is precompact in $H^s_0(\T)$.
Then, there exists a subsequence also denoted by $\{u_n\}_{n = 1}^\infty$
such that $u_n(0)$ converges to $u_\infty(0)$ in  $H^s_0(\T)$.
By the continuity  of the solution map $u(0) \in H^s_0(\T) \mapsto u \in X$,
 $u_n$ converges to $u_\infty : = \Phi (u_\infty(0))$ in $X$.
By weak convergence in $L^2_0(\T)$ of (a further subsequence of ) $\{u_n\}_{n = 1}^\infty$,
it is easy to see that  $u_\infty \in \wt K_r$.
This shows that $\wt K_r$ is compact in $X$.

By the continuity of the solution map $\Phi :  H^s_0(\T) \to X$, we
see that $\Psi(F) = \Phi^{-1}(F)$ is closed (and open) in $H^s_0(\T)$
if $F$ is closed (and open, respectively) in $X$.
Hence,
 as a direct consequence of \eqref{Z3} and \eqref{Z4},
we have the following large deviation principle
 for $\{\mu_\eps\}_{\eps > 0}$
with the rate function $\wt I$ defined in \eqref{Z5}.
We have
\begin{equation}
\limsup_{\eps \to 0} {\eps^2} \log \mu_\eps (F) \leq - \inf_{u \in F} \wt I(u)
\label{Z5a}
\end{equation}

\noi
for any closed set $F \subset X$
and 	
\begin{equation}
\liminf_{\eps \to 0} {\eps^2} \log \mu_\eps (G) \geq - \inf_{u \in G} \wt I(u)
\label{Z5b}
\end{equation}

\noi
for any open set $G \subset X$.
Let $v$ be  the unique solution to KdV
such that $v(0) = v_0$.
Since $\wt{K}_0 := \{ u \in X: \wt I(u) = 0\}$
consists of a single function $v$,
it follows again from Remark 2.3 in \cite{V2} that
 $\mu_\eps$ converges weakly to $\dl_v$.

\begin{remark}\rm
Due to lack of a good well-posedness theory below $L^2(\T)$,
the large deviation principle for
the Wick ordered cubic NLS \eqref{NLS1} holds  only in some mild sense.
Fix  $v_0 \in H^\al(\T)$ with $\al >\frac 5{12}$
and let $\phi(\omega)$ as in \eqref{FWb}.
With $s < \al -\frac{1}{2}$,
define a rate function $I: H^s(\T) \to [0, \infty]$
by
\begin{equation}
I(f) = \frac{1}{2} \|f - v_0\|_{H^\al(\T)}^2.
\label{Z6}
\end{equation}
	
\noi
Then, with $\rho_\eps = P\circ (u_0^\eps)^{-1}$,
the large deviation principle, i.e.~\eqref{Z3} and \eqref{Z4},
holds for $\{\rho_\eps\}_{\eps > 0}$
 as before.

However,
 the large deviation principle
for the probability measures $\mu_\eps$ induced by $\omega\mapsto u^\eps(\omega)$
holds only in a weak sense.
Let $X$ denote the collection of
all (known) solutions to \eqref{NLS1}
with initial data in $H^s(\T)$.
Due to lack of
a good well-posedness theory below $L^2(\T)$,
we do not know if $X$ is a Polish space and thus
we can only draw a weak conclusion for $\mu_\eps$.
With the previous notations,
we can define a rate function
$\wt I(u) := I(\Psi(u)) = I(u(0))$, where $I$ is as in \eqref{Z6}.
Then, we have the following `weak' large deviation principle
 for $\{\mu_\eps\}_{\eps > 0}$.
Namely, \eqref{Z5a} holds
if $\Psi(F)$ is closed in $H^s(\T)$,
while \eqref{Z5b} holds
if $\Psi(G)$ is open in $H^s(\T)$.
Lastly,
since $\wt{K}_0 := \{ u \in X: \wt I(u) = 0\}$
consists of a single function $v$,
where $v$ is the unique solution to \eqref{NLS1} such that $v(0) = v_0$,
we would like to conclude  that $\mu_\eps$ converges weakly to $\dl_v$.
Such weak convergence, however,
does not follow at this point
 due to lack of continuous dependence
on $H^s(\T)$
of the solution map to \eqref{NLS1}
constructed in \cite{CO} and Theorem \ref{THM:NLS} above.

\end{remark}

\section{Invariant set of full measure for almost sure global existence
of the Wick ordered cubic NLS}
\label{SEC:NLS}

In this section, we present the proof of Theorem \ref{THM:NLS}.
Let $\Phi(t): u(0) \mapsto u(t)$ be the solution map of \eqref{NLS1},
sending an initial condition $u(0)$ to the solution $u(t)$ at time $t$,
and $S(t) = e^{-i t \dx^2}$ be the linear propagator for \eqref{NLS1}.
In the following, we fix $\al \in (\frac{5}{12}, \frac{1}{2}]$
and $s = \al - \frac{1}{2}- \eps< 0$, $\eps > 0$.
The almost sure global result in \cite{CO}
states that
there exists a set $\Si_0 = \Si_0(\al)$
with
\begin{equation}
\rho_\al (\Si_0) = 1
\label{CO0b}
\end{equation}

\noi
 such that
if  $\phi = \phi(\omega)  \in \Si_0$,
the corresponding solution $u(t) = \Phi(t) \phi$ exists globally.
Here, $\rho_\al$ is as in \eqref{Gauss1} and $\phi \in \Si_0$ can be represented by \eqref{FWb}
almost surely.
Note that $\phi$ is almost surely in $H^s(\T) \setminus H^{\al -\frac{1}{2}}(\T)$.
In particular, it is not in $L^2(\T)$ almost surely for $\al < \frac{1}{2}$.
In establishing
this result,
we exhibited
nonlinear smoothing under randomization of the initial data,
i.e.~
if  $\phi   \in \Si_0$,
then, although the linear solution $S(t) \phi$ is not in $L^2(\T)$ for any $t\in \R$ almost surely,
the nonlinear part $v(t) := \Phi(t) \phi - S(t)\phi$ of the solution
is in $ L^2(\T)$ for each $t \in \R$.
Once we restrict our attention to the local-in-time setting,
we know more properties about this flow.
We summarize the local-in-time properties of the flow. See \cite[Theorem 1]{CO}.

\begin{proposition}[Summary of the local result in \cite{CO}] \label{PROP:CO}
 Fix $\dl \ll 1$.
Then, there exists $\Omega_\dl \in \mathcal{F}$
with the following properties.

\noi
\textup{(i)} The complemental measure of $\Omega_\dl$ is small.
More precisely, we have
\begin{equation}
P(\Omega_\dl^c) < e^{-\frac{1}{\dl^c}}.
\label{CO0a}
\end{equation}

\noi
\textup{(ii)}
For each $\omega \in \Omega_\dl$, there exists a unique solution $u$ to \eqref{NLS1}
in
\[ S(t) \phi(\omega) + C([-\dl, \dl]; L^2(\T))\subset C([-\dl, \dl]; H^s(\T))\]

\noi
with initial condition $u|_{t = 0} = \phi(\omega)$, where $\phi(\omega)$ is given by \eqref{FWb}.

\noi
\textup{(iii)}
Let $\omega \in \Omega_\dl$
and
$u(t) = \Phi(t) \phi (\omega) = S(t) \phi(\omega) + v(t)$ be the solution to \eqref{NLS1}
constructed in \textup{(ii)}.
Then, there exist $C$ and $ \theta > 0$ such that
\begin{equation} \label{CO1}
 \|v\|_{X^{0, \frac{1}{2}+}_{[-\dl', \dl']}} \le C (\dl')^\theta
\end{equation}

\noi
for $\dl' \leq \dl$.
Here, $X^{0, \frac{1}{2}+}_{[-\dl', \dl']}$ denotes the local-in-time version
of the $X^{\s, b}$-space restricted onto the time interval $[-\dl', \dl']$ with $\s = 0$ and $b = \frac{1}{2}+$.
In particular, we have
\begin{equation} \label{CO2}
\sup_{t\in [-\dl', \dl']}  \|v(t)\|_{L^2(\T)} \le C' (\dl')^\theta
\end{equation}

\noi
for some $C' > 0$.

\noi
\textup{(iv)} Let $w_0 \in L^2(\T)$ with $\|w_0\|_{L^2} = m$.
Then, there exists positive $\dl' = \dl'(m, \dl )< \dl$
such that, for each $\omega \in \Omega_\dl$,  there exists
a unique solution
$u\in C([t_*-\dl', t_*+\dl']; H^s(\T))$ to \eqref{NLS1}
with initial condition $u|_{t = t_*} = S(t_*) \phi(\omega) + w_0$
as long as $[t_*-\dl', t_*+\dl'] \subset [-\dl, \dl]$.
Here, uniqueness holds in the ball centered at
$S(\cdot - t_*) \big(S(t_*)\phi+w_0\big)$
of radius 1 in $X^{0, \frac{1}{2}+}_{[t_*-\dl', t_*+\dl']}$.

\end{proposition}

\noi
While (i) and (ii) of Proposition \ref{PROP:CO}
are exactly as in \cite[Theorem 1]{CO},
  (iii) and  (iv)  follow directly from (a modification of) the proof of \cite[Theorem 1]{CO}.
In particular, (iv) holds since
the required probabilistic estimates for the local argument in \cite{CO}
hold uniformly for any subinterval
$[t_*-\dl', t_*+\dl'] \subset [-\dl, \dl]$ if $\omega \in \Omega_\dl$.
We point out that we do not know if an analogue of (iv) holds
for the global-in-time setting.
Namely, given $w_0 \in L^2$,
(a small modification of) the proof of \cite[Theorem 2]{CO} does not yield
almost sure global existence for \eqref{NLS1}
with initial data $u|_{t = 0} = \phi(\omega) + w_0$,
where $\phi$ is as in \eqref{FWb}.

\medskip

In the following, we construct a set $\Si$ of full measure,
which is invariant under the \eqref{NLS1}-dynamics.
Moreover, our construction yields
an enhanced uniqueness statement (see Theorem \ref{THM:NLS} (ii)).

\smallskip

\noi
$\bullet$ {\bf Step 1:}
First, we use the invariance of $\rho_\al$ under the {\it linear} flow
and construct a set $\wt\Si$ of full measure
such that the linear solutions with initial data in $\wt \Si$
have some desired property.

For small $\dl > 0$, let $\Si_\dl = \phi(\Omega_{\dl})$,
where $\phi: \Omega \to H^s(\T)$ is the map given by \eqref{FWb}
and $\Omega_{\dl}$ is as in Proposition \ref{PROP:CO}.
Note that solutions with initial data in $\Si_{\dl}$ satisfy a good (local-in-time) uniqueness
property, coming from the local argument in \cite{CO}.
Letting
\[ \wt \Si_n := \Si_0 \cap \Si_\frac{1}{n},\]

\noi
we have $\rho_\al(\wt \Si_n^c) < e^{-n^c}$ for $n \geq N$.
Here, $N$ is a sufficiently large integer
such that Proposition \ref{PROP:CO} holds for all positive $\dl < N^{-1}\ll 1$.
Then, define $\ft \Si_n$ by
\[ \ft \Si_n := \bigcap_{k = 0}^n S\Big(-\frac{k}{n}\Big)\big(\wt \Si_n\big).\]

\noi
Since $\rho_\al$ is invariant under the  linear flow,
we have
\begin{equation}\label{Q1}
\rho_\al(\ft \Si_n^c) < (n+1) e^{-n^c}.
\end{equation}

\noi
Next, define $\wt \Si_{[0, 1]}$ by
\[ \wt \Si_{[0, 1]} := \bigcup_{n = N}^\infty \ft \Si_n.\]

\noi
Note that $\rho_\al(\wt \Si_{[0, 1]}) = 1$,
since
 \eqref{Q1} yields
\[ \rho_\al\big((\wt \Si_{[0, 1]})^c\big)
\leq \inf_{n \geq N} (n+1) e^{-n^c} = 0.\]

\noi
Finally, define $\wt \Si$ by
\[ \wt \Si := \bigcap_{j \in \Z} S(-j) \wt \Si_{[0, 1]}.\]

\noi
Then, by the invariance of $\rho_\al$ under the linear flow, we have
$ \rho_\al(\wt \Si) = 1$.

We claim that  if $\phi \in \wt \Si$, then given $t_* \in [j, j+1] \subset \R$,
 the conclusion
  of Proposition \ref{PROP:CO} (iv) holds
with $\dl = \frac{1}{n}$ for some $n = n(\phi,  j) \geq N$.
More precisely,
by writing $t_* = j + \tau$ with $\tau \in [0, 1]$,
we have $S(t_*) \phi = S(\tau) \psi$
for some $\psi = S(j) \phi \in \wt \Si_{[0, 1]}$,
i.e.~$\psi \in \ft{\Si}_n$ for some $n = n(\phi, j)$.
By further writing
 $\tau = \frac{k}{n} + \xi$ with $|\xi | < \frac{1}{2n}$
for some $k \in \{0, \dots, n\}$,
 we have
 $\varphi = S(\frac{k}{n}) \psi \in \wt \Si_n \subset \Si_{\frac 1n}$.
Then,
by Proposition \ref{PROP:CO} (iv),
given $w_0 \in L^2(\T)$ with $\|w_0\|_{L^2} = m$,
there exists $\dl' = \dl'(m, \dl) >0$ (with $\dl = \frac{1}{n}$)
such that  a solution $u$ to \eqref{NLS1}
with $u|_{t = \xi} = S(\xi) \varphi + w_0$
exists on $[\xi-\dl', \xi+\dl']$
and
is unique
in the ball centered at
$S(\cdot) \big(S(\xi)\varphi+w_0\big)$
of radius 1
 in $X^{0, \frac{1}{2}+}_{[\xi-\dl', \xi+\dl']}$.
Here, we assumed that $\dl' < \frac{\dl}{2} = \frac{1}{2n}$
such that
the subinterval $[\xi-\dl', \xi+\dl']$ lies in $[-\dl, \dl] = [-\frac 1n, \frac 1n]$.
Finally, note that
\[ S(\xi)\varphi = S(\tau) \psi
= S(t_*)\phi.\]

\noi
Therefore,
it follows from the  discussion above
that
there exists a unique solution $u$ to \eqref{NLS1}
on $[t_* - \dl', t_*+\dl']$
with $u|_{t = t_*} = S(t_*) \phi + w_0$,
where uniqueness holds
in the ball centered at
$S(\cdot) \big(S(t_*)\phi+w_0\big)$
of radius 1
 in $X^{0, \frac{1}{2}+}_{[t_*-\dl', t_*+\dl']}$.

\medskip

\noi
$\bullet$ {\bf Step 2:}
Next,  we construct the desired set $\Si$.
Define $\Si$ by
\[ \Si = \bigcup_{t \in \R} \Phi(-t) (\Si_0 \cap \wt \Si).\]

\noi
Recall that $\rho_\al$ is defined on the completion of
the Borel $\s$-algebra on $H^s(\T)$.  See Remark \ref{REM:meas} below.
Since $\Si \supset \Si_0 \cap \wt{\Si}$ with  $\rho_\al(\Si_0 \cap \wt{\Si}) = 1$,
it follows that  $\Si$ is measurable
and $\rho_\al(\Si) = 1$.
By definition, the set $\Si$ is invariant under the \eqref{NLS1}-dynamics.
Moreover, if
 $\phi \in \Si$, we have a global solution
\begin{equation}
u(t) =  \Phi(t) \phi = S(t) \phi + v(t),
 \label{CO3}
 \end{equation}

 \noi
where $v(t) \in L^2(\T)$ for each $t\in\R$.

Given $\phi \in \Si$, we have $\phi = \Phi(-t_0) \psi$
for some $t_0 \in \R$ and $\psi \in \Si_0 \cap \wt \Si$.
In particular,
given $t_* \in [j, j+1] \subset \R$,
we have
\[ u(t_*) =\Phi(t_*) \phi =  \Phi(t_* - t_0) \psi = S(t_* - t_0) \psi + w(t_*)\]

\noi
for some $w(t_*) \in L^2(\T)$.
Then,
by Step 1,
there exist
$\dl = \dl(\phi, t_*) = \dl(\phi, j) > 0$
and
 $\dl' = \dl'(t_*) = \dl'( \|w(t_*)\|_{L^2},  \dl) > 0$ such that
uniqueness  holds in
the ball of radius 1 centered at
$S(\cdot) u(t_*)$
in  $X^{0, \frac{1}{2}+}_{[t_*-\dl', t_*+\dl' ]}$.
Lastly, given a finite time interval $I \subset [J_1, J_2] \subset \R$,
we have $\inf_{t_* \in I} \dl' (t_*)>0$
since $\sup_{t_* \in I}\|w(t_*)\|_{L^2} < \infty$
and  $\inf_{t_*\in I} \dl(t_*) =  \inf \{ \dl(\phi, j): j = J_1,J_1+1,  \dots, J_2\}>0$.
This completes the proof of Theorem \ref{THM:NLS}.

\begin{remark}\label{REM:meas}
\rm
The Gaussian measure $\rho_\al$
is the induced probability measure
under the map $\phi: \omega \in \Omega \mapsto \phi^\omega \in H^s(\T)$, $s < \al -\frac 12$,  defined in \eqref{FWb}.
In the following, we directly show that $\rho_\al$ is defined on the Borel $\s$-algebra in $H^s(\T)$.
Let $\phi_N$ be the Fourier truncation of $\phi$ given by
\begin{equation} \label{phi_N}
\phi_N(x; \omega) = \sum_{|n|\leq N} \frac{g_n(\omega)}{\jb{|n|^{\al}}}\,  e^{inx}.
\end{equation}

\noi
Then,
the set $A_{N, r} = \{ \omega \in \Omega: \|\phi_N(\omega)\|_{H^s} \leq  r\}$
is clearly measurable for each $N\in \mathbb{N}$ and for any $r\geq 0$.
With
 $A_r = \{ \omega \in \Omega: \|\phi(\omega)\|_{H^s} \leq  r\}$,
 we have $A_r = \bigcap_{N\in \mathbb{N}} A_{N, r}$
 and hence $A_r$ is also measurable.
Let  $ B_r(v)$ be the open ball of radius $r$ centered at $v\in H^s(\T)$.
Then, by writing $ B_r(v) = \bigcup_{n = 1}^\infty \cj B_{r - \frac 1n}(v)$,
we see that
$\phi^{-1}( B_r(v))$ is measurable.
Since $H^s(\T)$ is separable,
any open set can be written as a countable union
of open balls in $H^s(\T)$.
Hence, we conclude that $\rho_\al$ is defined on the Borel $\s$-algebra in $H^s(\T)$
(and on its completion).
\end{remark}


\smallskip
\noi
{\bf Acknowledgments.}
T.O.~
would like to thank Prof.~D.~Calegari
for asking a question related to Theorem \ref{THM:NLS}
during his presentation at the University of Cambridge.
He is grateful to Prof.~J.~Colliander for interesting discussions
and his encouragement.
He also thanks G. Richards for a useful discussion on
the variational characterization of Gibbs measures.
The authors
are grateful to Prof.~Y.~Tsutsumi
for asking a question related to the large deviation principle.

\end{document}